\documentclass[11pt]{article}
\usepackage{amsmath}
\usepackage{amsfonts}
\usepackage[bookmarksnumbered=true]{hyperref}

\newtheorem{theorem}{Theorem}[section]
\newtheorem{lemma}[theorem]{Lemma}
\newtheorem{definition}[theorem]{Definition}
\newtheorem{proposition}[theorem]{Proposition}
\newtheorem{corollary}[theorem]{Corollary}
\newtheorem{remark}[theorem]{Remark}
\def\dfrac{\displaystyle\frac}
\def\dint{\displaystyle\int}

\begin{document}

\title{A Kneser-type theorem for backward doubly stochastic differential
equations\thanks{This work is supported by National Natural Science Foundation of China Grant
10771122, Natural Science Foundation of Shandong Province of China
Grant Y2006A08 and National Basic Research Program of China (973
Program, No.2007CB814900)}}
\author{Yufeng Shi$^{\rm a}$\thanks{Corresponding author, E-mail: yfshi@sdu.edu.cn}\ and
 Qingfeng Zhu$^{\rm b}$\\
{\small$^{\rm a}$School of Mathematics, Shandong University, Jinan 250100, China}\\
{\small $^{\rm b}$ School of Statistics and Mathematics, Shandong
University of Finance},\\
{\small Jinan 250014, China}}
\maketitle

\begin{abstract}A class of backward doubly stochastic differential equations (BDSDEs
in short) with continuous coefficients is studied. We give the
comparison theorems, the existence of the maximal solution and the
structure of solutions for BDSDEs with continuous coefficients. A
Kneser-type theorem for BDSDEs is obtained. We show that there is
either unique or uncountable solutions for this kind of BDSDEs.\\
\indent{\it keywords:} Backward doubly stochastic differential equations, 
comparison theorem, maximal solution, Kneser-type theorem
\end{abstract}

\section{Introduction}

Nonlinear backward stochastic differential equations (BSDEs in
short) have been independently introduced by Pardoux and Peng
\cite{PP1} and Duffie and Epstein \cite{DE}. The comprehensive
applications of BSDEs have motivated many efforts to establish the
existence and uniqueness of adapted solution under general
hypotheses on the coefficients. For instance, for the
one-dimensional case, Lepeltier and San Martin \cite{LM} proved the
existence of a solution to BSDEs under the assumption of continuous
coefficient. Recently, Jia and Peng \cite{JP} studied the structure of the
solutions to BSDEs with continuous coefficients.

The comparison theorem is an important and effective technique in
the theory of BSDEs. There are much works concerning the comparison
theorems of BSDEs. For instance, for applications to finance, El
Karoui, Peng and Quenez \cite{KPQ} have given the comparison theorem
 of BSDEs with Lipschitz coefficients. And then Liu
and Ren \cite{LR} have given some comparison theorems of BSDEs with continuous coefficients.

In 1923, Kneser \cite{Kn} proved that the cardinality of
the set of solutions for ordinary differential equations (ODEs in short)
with continuous coefficients is either one or of continuum.
Many interesting generalizations have followed. For example,
in 1956, Alexiewicz and Orlicz \cite{AO} got the same theorem for
a class of partial differential equations (PDEs in short).
In 2008, Jia \cite{J} proved Kneser-type theorem
for BSDEs with continuous coefficients. For more
information about Kneser-type theorems, one can refer
to \cite{K, Ke, S, V}.

A class of backward doubly stochastic differential equations (BDSDEs
in short) was introduced by Pardoux and Peng \cite{PP2} in 1994, in
order to provide a probabilistic interpretation for the solutions of
a class of semilinear stochastic partial differential equations
(SPDEs in short). They have proved the existence and uniqueness of
solution for BDSDEs under uniformly Lipschitz conditions. Since
then, many efforts have been made to relax the assumption on the
coefficients. For instance, Gu \cite{G} proved the existence and
uniqueness of solution for BDSDEs under the local Lipschitz
conditions. Shi, Gu and Liu \cite{SGL} have relaxed the Lipschitz
assumptions to linear growth conditions by virtue of their
comparison theorem of BDSDEs with Lipschitz conditions, and showed
the existence of the minimal solution of BDSDEs. Bally and Matoussi
\cite{BM} have given a probabilistic interpretation of the solutions
in Sobolev spaces for parabolic semilinear stochastic PDEs in terms
of BDSDEs. Zhang and Zhao \cite{ZZ} have proved the existence and
uniqueness of solution for BDSDEs on infinite horizon, and
described the stationary solutions of SPDEs by virtue of the
solutions of BDSDEs on infinite horizon.

Due to their important significance to SPDEs, it is necessary to
give intensive investigation to the theory of BDSDEs. The aim of
this paper is to study the structure of the solutions to BDSDEs with
continuous and linear growth conditions. We firstly generalize the
comparison theorems to the case where the coefficients are
continuous. As an application of the comparison theorems, we give
the existence of the maximal solution of BDSDEs with continuous
coefficient by means of our comparison theorems. Finally we will
show that there is either unique or uncountable solutions for this
kind of BDSDEs. In fact, our result shows the structure of those
solutions, that is, we obtain a Kneser-type theorem for BDSDEs.

The rest of the paper is organized as follows. In Section 2, we
present the main assumptions and some preliminary results. In
Section 3, we generalize the comparison theorem of BDSDEs in
\cite{SGL} to BDSDEs with continuous coefficients. In Section 4, we
prove the existence of the maximal solution of BDSDEs. Finally, in
Section 5, we discuss the structure of solutions of BDSDEs with
continuous coefficients and linear growth conditions.

%The title of your section 2
\section{Preliminaries}

Let $(\Omega,\mathcal {F},P)$ be a probability space, and $T>0$ be
an arbitrarily fixed constant throughout this paper. Let
$\{W_t;0\leq t\leq T\}$ and $\{B_t;0\leq t\leq T\}$ be two mutually
independent standard Brownian Motions with values in $\mathbf{R}^d$
and $\mathbf{R}^l$, respectively, defined on
$(\Omega,\mathcal{F},P)$. Let $\mathcal{N}$ denote the class of
$P$-null sets of $\mathcal{F}$. For each $t\in [0,T]$, we define
${\mathcal{F}}_t \doteq {\mathcal{F}}_t^W\vee
{\mathcal{F}}_{t,T}^B$, where for any process $\{\eta_t\}$,
${\mathcal{F}}_{s,t}^{\eta}=\sigma\{\eta_r-\eta_s;s\leq r\leq
t\}\vee {\mathcal{N}}$,
${\mathcal{F}}_t^{\eta}={\mathcal{F}}_{0,t}^{\eta}$. Note that the
collection $\{{\mathcal{F}}_t;t\in [0,T]\}$ is neither increasing
nor decreasing, so it does not constitute a filtration.
All the equalities and inequalities mentioned in this paper are
in the sense of ${\rm d}t\times {\rm d}P$ almost surely on $\left[
0,T\right] \times \Omega .$

We use the usual inner
product $\left\langle \cdot ,\cdot \right\rangle $ and Euclidean norm $
\left| \cdot \right| $ in $\mathbf{R}^k, \mathbf{R}^{k\times l}$
and $\mathbf{R}^{k\times
d}.$ All the equalities and inequalities mentioned in this paper are
in the sense of ${\rm d}t\times {\rm d}P$ almost surely on $\left[
0,T\right] \times \Omega .$

For any $k\in N$, let $M^2(0,T;\mathbf{R}^k)$ denote the set of
(classes of ${\rm d}P\otimes {\rm d}t$ a.e. equal) $k$-dimensional
jointly measurable stochastic processes $\{\varphi_t;t\in [0,T]\}$
which satisfy:

\noindent(i)\ $\| \varphi \|_{M^2}^2 :=
\mathbf{E}\int_0^T|\varphi_t|^2{\rm d}t<\infty$;

\noindent(ii)\ $\varphi_t$ is ${\mathcal{F}}_t$-measurable, for any
$t\in
[0,T]$.\\
Similarly, we denote by $S^2([0,T];\mathbf{R}^k)$ the set of
$k$-dimensional continuous stochastic processes $\{\varphi_t;t\in
[0,T]\}$ which satisfy:

\noindent(iii)\ $\| \varphi \|_{S^2}^2:= \mathbf{E}(\sup_{0\leq
t\leq T}|\varphi_t|^2)<\infty$;

\noindent(iv)\ $ \varphi_t$ is ${\mathcal{F}}_t$-measurable, for any
$t\in[0,T]$.\\
For any $t\in[0,T]$, denote by
$L^2(\Omega,{\mathcal{F}}_t,P;\mathbf{R}^k)$ the set of
 $k$-dimensional random variables $\xi$, which satisfy:

\noindent(v)\ ${\bf E}|\xi|^2 <\infty$;

\noindent(vi)\ $\xi$ is ${\mathcal{F}}_t$-measurable.

Consider the following BDSDE:
\begin{equation}\label{1}
 Y_t=\xi+\int_t^Tf(s,Y_s,Z_s){\rm d}s+\int_t^Tg(s,Y_s,Z_s){\rm
d}B_s-\int_t^TZ_s{\rm d}W_s,\quad 0\leq t\leq T.
\end{equation}
We assume that

(H1)\ $f:\Omega\times [0,T]\times \mathbf{R}^k\times \mathbf{R}^{k\times
d}\rightarrow \mathbf{R}^k$ is jointly measurable and such that
\begin{align*}f(\cdot,y,z)\in M^2(0,T;\mathbf{R}^k),
\ \forall (y,z)\in \mathbf{R}^k\times \mathbf{R}^{k\times d},
\end{align*}
and there exist a constant $C>0$ such that for any $t \in [0,T]$,
  $(y_i,z_i)\in \mathbf{R}^k\times \mathbf{R}^{k\times d}$, $i=1,2$,
\begin{align*}|f(t,y_1,z_1)-f(t,y_2,z_2)|^2 \leq C(|y_1-y_2|^2+|z_1-z_2|^2).\end{align*}

(H2)\ $g:\Omega\times [0,T]\times \mathbf{R}^k\times \mathbf{R}^{k\times d}\rightarrow
\mathbf{R}^{k\times l}$ is jointly measurable and such that
\begin{align*}g(\cdot,y,z)\in M^2(0,T;\mathbf{R}^k),\
\forall (y,z)\in \mathbf{R}^k\times \mathbf{R}^{k\times d},\end{align*}
and there exist constants $C>0$ and $0<\alpha<1$
such that for any $t \in [0,T]$, $(y_i,z_i)\in \mathbf{R}^k\times \mathbf{R}^{k\times d}$, $i=1,2$,
\begin{align*}|g(t,y_1,z_1)-g(t,y_2,z_2)|^2\leq C|y_1-y_2|^2+\alpha |z_1-z_2|^2.\end{align*}

Note that the integral with respect to $\{B_t\}$ is a ``backward
It\^o integral", in which the integrand takes values
at the right end points of the subintervals in the
Riemann type sum (for details refer to \cite{PP2}), and the integral with respect to $\{W_t\}$ is a
standard forward It\^o integral. These two types of integrals are
particular cases of the It\^o-Sokorohod integral (see \cite{NP} for
details).

\begin{definition} A pair of
processes $(Y,Z):\Omega \times [0,T]\rightarrow \mathbf{R}^k\times
\mathbf{R}^{k\times d}$ is called a solution of BDSDE (1), if
$(Y,Z)\in S^2([0,T];\mathbf{R}^k)\times M^2(0,T;\mathbf{R}^{k\times
d})$ and satisfies BDSDE (1).
\end{definition}

\begin{proposition} Under assumptions {\rm (H1)} and {\rm(H2)},
 if $\xi\in L^2(\Omega,{\mathcal{F}}_T,P;\mathbf{R}^k)$, then BDSDE {\rm (1)} has
a unique solution
\begin{align*}
(Y,Z)\in S^2([0,T];\mathbf{R}^k)\times M^2(0,T;\mathbf{R}^{k\times
d}).
\end{align*}
\end{proposition}
\noindent This proposition was derived in \cite{PP2}.

%The title of your section 2
\section{Comparison theorem of BDSDEs with continuous coefficients}

In this paper, we only consider one-dimensional BDSDEs, i.e., $k=1$.
Assume

(H3)\ for fixed $\omega$ and $t$, $f(\omega, t, \cdot,\cdot)$ is
continuous;

(H4)\ linear growth: $\forall ( \omega, t, y, z)\in \Omega \times
[0,T]\times \mathbf{R}\times \mathbf{R}^d$, there exists $0<K<
\infty$, such that
\begin{align*}
|f(\omega,t,y,z)|\leq|f(\omega,t,0,0)|+K|y|+K|z|,
\end{align*}
 and  \begin{align*}
\mathbf{E}\int_0^T|f(\omega,t,0,0)|^2{\rm dt}<\infty.
\end{align*}

\vspace{1mm} We consider the following BDSDEs:\hspace{0.1cm}($0\leq
t\leq T$)
\begin{equation}
 Y^1_t=\xi^1+\int_t^Tf^1\left(s,Y^1_s,Z^1_s\right)
{\rm d}s+\int_t^Tg\left(s,Y^1_s,Z^1_s\right){\rm
d}B_s-\int_t^TZ^1_s{\rm d}W_s,
\end{equation}
\begin{equation}
Y^2_t=\xi^2+\int_t^Tf^2\left(s,Y^2_s,Z^2_s\right) {\rm
d}s+\int_t^Tg\left(s,Y^2_s,Z^2_s\right){\rm d}B_s-\int_t^TZ^2_s{\rm
d}W_s.
\end{equation}
where for $i=1,2$, $\xi^i \in
L^2(\Omega,{\mathcal{F}}_T,P;\mathbf{R})$ and $f^i:\Omega\times
[0,T]\times \mathbf{R}\times \mathbf{R}^{d}\rightarrow \mathbf{R}$
satisfy

(H5) $\xi^1\geq\xi^2$,\ a.s.; \quad$f^1(t,y,z)\geq f^2(t,y,z),\
a.s., \ \forall\ (t,y,z)\in [0,T]\times \mathbf{R}\times
\mathbf{R}^d.$

The comparison theorem was established by Shi, Gu and Liu \cite{SGL}
for one-dimensional BDSDEs, where both the coefficients
$f^1$ and $f^2$ satisfy Lipschitz conditions. In this section, we
firstly generalize the comparison theorem to the case where one of
the coefficients $f^1$ and $f^2$ is only continuous, another is
Lipschitz continuous.

\begin{theorem}  \label{result1}
Assume BDSDE {\rm (2)} satisfies {\rm (H1)} and {\rm (H2)}, and
BDSDE {\rm (3)} satisfies {\rm (H2)}--{\rm (H4)}. Let $(Y^1,Z^1)$
and $(Y^2,Z^2)$ be the solutions of BDSDEs {\rm (2)} and {\rm (3)},
respectively. If {\rm (H5)} holds, then
\begin{align*}P\{Y_t^1\geq Y_t^2,\hspace{0.2cm} \mbox{for all }t \in
[0,T]\}=1.\end{align*}
\end{theorem}

{\bf Proof.}\
 It is easy to see that
$(Y^1_t- Y^2_t, Z^1_t- Z^2_t)$ satisfies the following BDSDE
\begin{equation*}\label{Multi}
 \begin{split}
    Y_t^1-Y_t^2 =&  \xi^1-\xi^2+\int_t^T[f^1\left(s,Y^1_s,Z^1_s\right)
    -f^2\left(s,Y^2_s,Z^2_s\right)]{\rm d}s\\
    &  +\int_t^T[g\left(s,Y^1_s,Z^1_s\right)-g\left(s,Y^2_s,Z^2_s\right)]{\rm d}B_s\\
    &  -\int_t^T(Z_s^1-Z_s^2){\rm d}W_s,\quad 0\leq t\leq T.
 \end{split}
\end{equation*}
Applying It\^{o}-Tanaka's formula (cf. \cite{SGL}) to $|(Y_s^1-Y_s^2)^-|^2$, we get
\begin{equation*}\label{Multi}
 \begin{split}
    & |(Y_t^1-Y_t^2)^{-}|^{2}\\
    =&  |(\xi^1-\xi^2)^-|^2
    -2\int_t^T(Y_s^1-Y_s^2)^-[f^1\left(s,Y^1_s,Z^1_s\right)
    -f^2\left(s,Y^2_s,Z^2_s\right)]{\rm d}s\\
    &  -2\int_t^T(Y_s^1-Y_s^2)^-
    [g\left(s,Y^1_s,Z^1_s\right)-g\left(s,Y^2_s,Z^2_s\right)]{\rm d}B_s\\
    &  +\int_t^T{\bf 1}_{\{Y_s^1\leq Y_s^2\}}
    |g\left(s,Y^1_s,Z^1_s\right)-g\left(s,Y^2_s,Z^2_s\right)|^2{\rm d}s\\
   & +2\int_t^T(Y_s^1-Y_s^2)^-(Z_s^1-Z_s^2){\rm d}W_s-\int_t^T{\bf
1}_{\{Y_s^1\leq Y_s^2\}}|(Z_s^1-Z_s^2)|^2{\rm d}s.
 \end{split}
\end{equation*}
From (H5), we have $\xi^1-\xi^2\geq 0$, so
\begin{equation}
\nonumber{\bf E}|(\xi^1-\xi^2)^-|^2=0.
\end{equation}
Since $(Y^1,Z^1)$ and $(Y^2,Z^2)$ are in $S^2(0,T;{\bf R})\times
M^2(0,T;{\bf R}^{d})$, by virtue of Lemma 1.3 of \cite{PP2}, it follows that
\begin{equation}
\nonumber{\bf E}\int_t^T(Y^1_s-Y^2_s)^-(Z^1_s-Z^2_s){\rm d}W_s=0,
\end{equation}
\begin{equation}
\nonumber{\bf
E}\int_t^T(Y^1_s-Y^2_s)^-[g\left(s,Y^1_s,Z^1_s\right)-g\left(s,Y^2_s,Z^2_s\right)]{\rm
d}B_s=0.
\end{equation}
Let
\begin{align*}
\Delta  =  -2\int_t^T(Y^1_s-Y^2_s)^-[f^1\left(s,Y^1_s,Z^1_s\right)
-f^2\left(s,Y^2_s,Z^2_s\right)]{\rm d}s = \Delta_1+\Delta_2,
\end{align*}
where
\begin{align*}
\Delta_1&=  -2\int_t^T(Y^1_s-Y^2_s)^-[f^1\left(s,Y^1_s,Z^1_s\right)
-f^1\left(s,Y^2_s,Z^2_s\right)]{\rm d}s,\\
\Delta_2&=  -2\int_t^T(Y^1_s-Y^2_s)^-[f^1\left(s,Y^2_s,Z^2_s\right)
-f^2\left(s,Y^2_s,Z^2_s\right)]{\rm d}s\leq 0.
\end{align*}
From (H1) and Young's inequality, it follows that
\begin{align*}
\Delta &\leq \Delta_1 \leq
2c\int_t^T(Y^1_s-Y^2_s)^-\left(|Y^1_s-Y^2_s|+|Z^1_s-Z^2_s|\right){\rm d}s \\
&\leq
\left(2c+\frac{c^2}{1-\alpha}\right)\int_t^T|(Y^1_s-Y^2_s)^-|^2{\rm
d}s +\left(1-\alpha\right)\int_t^T{\bf1}_{\{Y^1_s\leq
Y^2_s\}}|Z^1_s-Z^2_s|^2{\rm d}s,
\end{align*}
where $c>0$ only depends on the Lipschitz constant $C$ in (H1). By
the assumption (H2), we deduce
\begin{align*}
&\quad \int_t^T{\bf 1}_{\{Y^1_s\leq Y^2_s
\}}|g\left(s,Y^1_s,Z^1_s\right)
-g\left(s,Y^2_s,Z^2_s\right)|^2{\rm d}s\\
&\leq \int_t^T{\bf 1}_{\{Y^1_s\leq Y^2_s\}}
\left[C|Y^1_s- Y^2_s|^2+\alpha|Z^1_s- Z^2_s|^2\right]{\rm d}s\\
& = C\int_t^T|(Y^1_s- Y^2_s)^-|^2{\rm d}s +\alpha\int_t^T{\bf
1}_{\{Y^1_s\leq Y^2_s\}}|Z^1_s- Z^2_s|^2{\rm d}s.
\end{align*}
We get
\begin{equation*}{\bf E}|(Y^1_t- Y^2_t)^-|^2\leq
\left(C+2c+\frac{c^2}{1-\alpha}\right) \int_t^T{\bf E}|(Y^1_s-
Y^2_s)^-|^2{\rm d}s.
\end{equation*}
By Gronwall's inequality, it follows that
\begin{equation*}
{\bf E}|(Y^1_t- Y^2_t)^-|^2=0, \forall\ t \in [0,T].
\end{equation*}
That is, $P\{Y_t^1\geq Y_t^2$, for all $t \in [0,T]\}=1.$
\quad$\Box$

\begin{remark}
If BDSDE {\rm (2)} satisfies {\rm (H2)}--{\rm (H4)}, and BDSDE {\rm
(3)} satisfies {\rm (H1)} and {\rm (H2)}, similarly to Theorem 3.1,
the comparison result is still true.
\end{remark}

In next section, we will prove the maximal solution of BDSDE with
continuous coefficients as an application of Theorem 3.1. Next, let
us generalize the comparison theorem to the case where the
coefficients are only continuous.

\begin{theorem}  \label{result2}
Assume BDSDEs {\rm (2)} and {\rm (3)} satisfy {\rm (H2)-(H4)},
respectively. Let $(\underline Y^1,\underline Z^1)$ and $(\underline
Y^2,\underline Z^2)$ be the minimal solutions of BDSDEs {\rm (2)}
and {\rm (3)}, respectively. If {\rm (H5)} holds, then
\begin{equation*}
P\{\underline Y_t^1\geq \underline Y_t^2,\hspace{0.2cm}
\mbox{\rm for all }t \in [0,T]\}=1.
\end{equation*}

\noindent Moreover, let $(\overline Y^1,\overline Z^1)$ and
$(\overline Y^{\hspace{1mm}2},\overline Z^{\hspace{1mm}2})$ are the
maximal solutions of BDSDEs {\rm (2)} and {\rm (3)}, respectively.
If {\rm (H5)} holds, then
\begin{equation*}
P\{\overline  Y_t^1\geq \overline
Y_t^{\hspace{1mm}2},\hspace{0.2cm} \mbox{\rm for all }t \in
[0,T]\}=1.
\end{equation*}
\end{theorem}

{\bf Proof.}\ First, for fixed $t$, we define, as in Lemma 1 of \cite{LM},
 the sequence $f^2_n(t,y, z)$ associated with $f^2$
$$f^2_n(t,y, z)= \inf_{(y', z')\in Q\times Q^{d}}\{f^2(t,y',z')+n(|y-y'|+|z-z'|)\},$$
where $Q$ is the rational number set. Then, for $n \geq K$, $f^2_n$ are
measurable and Lipschitz functions, and $ f^1\geq f^2\geq f^2_n$ .
Hence, we know that the following BDSDE has a unique solution
$(Y^2_n, Z^2_n)$
\begin{equation*}
 \begin{split}
Y^2_n(t)=& \xi^2+\int_t^Tf^2_n\left(s,Y^2_n(s),Z^2_n(s)\right) {\rm
d}s+\int_t^Tg\left(s,Y^2_n(s),Z^2_n(s)\right){\rm
d}B_s\\
&  -\int_t^TZ^2_n(s){\rm d}W_s,\ n \geq K.
 \end{split}
\end{equation*}
From Theorem 3.1, it follows that $Y^1\geq Y^2_n$ and $Y^2 \geq
Y^2_n$ a.s. for all $n\geq K$. However, $(Y^2_n, Z^2_n)$ converges
uniformly in $t$ to $(\underline Y^2, \underline Z^2)$ (cf. \cite{SGL}). Therefore, $Y^1 \geq \underline Y^2$ a.s., in
particular, $\underline Y^1 \geq \underline Y^2$ a.s..

\ Next, for fixed $t$,  we define the sequence $f^1_n(t,y, z)$ associated with $f^1$,
$$f^1_n(t,y, z)= \sup_{(y',
z')\in Q\times Q^{d}}\{f^1(t,y',z')-n(|y-y'|+|z-z'|)\}$$
then, by virtue of Lemma 4.2 in next section, for $n \geq K$, $f^1_n$ are
measurable and Lipschitz functions, and $f^1_n\geq f^1\geq f^2$.
Hence, we know that the following BDSDE has a unique solution
$(Y^1_n, Z^1_n)$
\begin{equation*}
 \begin{split}
Y^1_n(t)=&\xi^1+\int_t^Tf^1_n\left(s,Y^1_n(s),Z^1_n(s)\right) {\rm
d}s+\int_t^Tg\left(s,Y^1_n(s),Z^1_n(s)\right){\rm
d}B_s\\
& -\int_t^TZ^1_n(s){\rm d}W_s,\ n \geq K.
 \end{split}
\end{equation*}
From Theorem 3.1, it follows that $Y^1_n \geq Y^1$ and $Y^1_n \geq
Y^2$ a.s. for all $n\geq K$. However, $(Y^1_n, Z^1_n)$ converges
uniformly in $t$ to $(\overline Y^1, \overline Z^1)$ (cf. Lemma
4.2). Therefore, $ \overline Y^1 \geq Y^{2} $ a.s., in particular,
$\overline Y^1 \geq \overline Y^{\hspace{1mm}2}$ a.s..
\quad$\Box$

\begin{remark} Suppose that the conditions of Theorem 3.2 hold.

(i) If (2) has a unique solution $(Y^1, Z^1)$, then for any
solutions $(Y^2, Z^2)$ of (3) we have
\begin{align*}
P\{ Y_t^1\geq Y_t^2,\hspace{0.2cm} \mbox{for all }t \in [0,T]\}=1.
\end{align*}
\indent(ii) If (3) has a unique solution $(Y^2, Z^2)$, then for any
solutions $(Y^1, Z^1)$ of (2) we have
\begin{align*}
P\{ Y_t^1\geq Y_t^2,\hspace{0.2cm} \mbox{for all }t \in [0,T]\}=1.
\end{align*}

\indent(iii) If the uniqueness of neither (2) nor (3) holds, we were
unable to derive a comparison result for any solutions of (2) and
(3).

\end{remark}
In particular, we easily see that the Lipschitz condition is a
special case of our proposed conditions. In other words, Theorem 3.1
and Theorem 3.2 generalize the comparison result in \cite{SGL}.

Finally, before closing this section, we shall give a comparison
result for the case of Remark 2(iii). For this purpose we need now
the following conditions:

(H6)\ $f(t, y, z)$  is uniformly Lipschitz continuous in $z$, i.e.,
$$|f(t, y, z_1)-f(t, y, z_2)|< C|z_1-z_2|$$
for all $t \in[0, T], y \in R, z_1, z_2 \in R^d$;

(H7)\ For fixed $t, z, f(t, \cdot, z)$ is equi-continuous in $y$,
i.e., for every $\varepsilon>0$ there exists a constant $\delta>0$
such that if $y_1, y_2\in R$, and $|y_1- y_2|<\delta$, then $|f(t,
y_1, z)-f(t, y_2, z)|\leq \varepsilon$ for any $t\in [0, T]$ and $z
\in R^d$;

(H8)\ $|\xi|\leq C_0$, a.s.

\noindent where $C$ and $C_0$ are given positive constants.

(H9)\ $f(t, y, z)$ is locally Lipschitz continuous in $y$, i.e., for
each $N>0$, there exists a constant $L_N$ such that
$$|f(t, y_1, z)-f(t, y_2, z)|\leq L_N |y_1-y_2|$$
for all $t \in[0, T]$, $z \in R^d$, $y_1$, $y_2 \in R$ with $|y_1|$,
$|y_2|\leq N$.

Moreover, by virtue of the results in \cite{G}, (1) is uniquely solvable under
the conditions: (H2), (H4), (H6), (H8) and (H9). In view of
Remark 2, if one of $f^1$ and $f^2$ satisfies (H4), (H6), (H8)
and (H9), and $g$ satisfies (H2), the conclusion of Theorem 3.2
still holds.

\begin{theorem}  \label{result3}
 Assume ${\rm (H2), (H4), (H8)}$ and one of $f^1$ and $f^2$ satisfies ${\rm
(H6)}$ and ${\rm (H7)}$ with

${\rm (H10)}$ $$ {\rm esssup}_{\{t\in[0,T], y\in R, z \in R^d\}}
\{f^1(t, y, z)-f^2(t, y, z)\}< \infty,$$ and let
$(Y^1, Z^1)$ and $(Y^2, Z^2)$ be the solutions of BDSDEs {\rm (2)}
and {\rm (3)}, respectively. Then {\rm (H5)} implies that
$$P\{Y^1(t)\geq Y^2(t), \mbox{\rm for all}\ t \in[0, T]\}= 1.$$
\end{theorem}
{\bf Proof.}\
The proof is divided into two steps. Without loss of generality, we
assume that $f^2$ satisfies (H6) and (H7). In the first step, we
additionally assume that $f^2(t, y, z)$ satisfies (H9) and shall
prove the result of Theorem 3.3 without the assumption (H7).

Under the assumptions (H2), (H4), (H6), (H8) and (H9), (1)
with $(\xi^2, f^2)$ has a unique solution by virtue of Theorem 2.1
of Gu \cite{G}. An application of Remark 2 yields the desired result.

The second step, it remains to remove the additional condition (H9).
Fix $N>0$ and choose a Lipschitz continuous mapping $b(t, y, z)$ in
$y$ and $z$ such that $$f^1(t, y, z)>b(t, y, z)>f^2(t, y, z)$$ for
$t \in[0,T]$, $y \in R$ with $|y|<N$, $z \in R^d$.

From (H10), there exists a constant $\bar \varepsilon >0$ such that
$${\rm esssup}_{\{t\in[0,T], y\in R, z \in R^d\}}
\{f^1(t, y, z)-f^2(t, y, z)\}=\bar \varepsilon.$$

To show the
existence of such a mapping $b(t, y, z)$, we write $\bar b(t, y, z)
:=f^2(t, y, z) +\displaystyle\frac{\bar\varepsilon}{2}$,

$$J(y)=\left\{
\begin{array}{cc}
k\exp(-(1-|y|)^{-1}),\ \ &\ \ \ \mbox{for}\  |y|<1,\\
0,\ \ &\ \ \ \mbox{otherwise},
\end{array}
\right.
$$
and $J_\delta (y)=\delta^{-1}J(y/\delta)$, where the constant $k$
satisfies $\int_{R}J(y)dy=1$. Let us smooth out $\bar b$ in $y$ to
obtain $b_\delta$, i.e., setting
$$b_\delta(t,y,z)=\dint_{R}\bar b(t,x,z)
J_\delta(y-x){\rm d}x,$$ where $\bar\varepsilon$ and $\delta$ are as in
(H7). Now any $b_{\delta'}$ with $\delta'\leq \delta$ is
our candidate.

Let $(Y, Z)$ denote a solution of (1) when $f$ is replaced by $b$. By
the above consideration, we observe $$Y^1(t)\geq Y(t)\geq
Y^2(t),a.s.\quad \forall t \in [0,T].$$ As a result of that $N>0$ is
arbitrary, the desired result is obtained.
\quad$\Box$

\begin{remark}
(i)\ Theorem 3.3 gives a comparison result for any solutions of (2)
and (3).

(ii)\ If one of (2) and (3) has a unique solution, Theorem 3.2
implies Theorem 3.3. Otherwise, Theorem 3.3 is not a special case of
Theorem 3.2.
\end{remark}

\section{The existence of the maximal solution of BDSDEs}

Under the conditions (H2)-(H4), Shi, Gu and Liu \cite{SGL} have
proved that BDSDE (1) has the minimal solution $(\underline Y,
\underline Z)$. Now, we will prove that BDSDE (1) has also the
maximal solution $(\overline Y, \overline Z)$.

\begin{theorem} \label{result4}
 Assume $f:\Omega\times [0,T]\times \mathbf{R} \times
\mathbf{R}^d\rightarrow \mathbf{R}$ and $g:\Omega \times [0,T]
\times \mathbf{R} \times \mathbf{R}^d\rightarrow \mathbf{R}^l$ are
jointly measurable functions and satisfy {\rm (H2)-(H4)}. Then, if
$\xi\in L^2(\Omega,{\mathcal{F}}_T,P;\mathbf{R})$, BDSDE {\rm (1)}
has a solution $(Y, Z)\in S^2([0,T];\mathbf{R})\times
M^2(0,T;\mathbf{R}^d)$. Moreover, there is a maximal solution
$(\overline{Y},\overline{Z})$ of BDSDE {\rm (1)} in the sense that,
for any other solution $(Y,Z)$ of BDSDE {\rm (1)}, we have
$\overline{Y}\geq Y$.
\end{theorem}

In order to prove Theorem 4.1, we need the following lemma which can be obtained by means of the similar arguments in
\cite{LM}, so we omit its proof.

\begin{lemma}
Let $f:\mathbf{R}\times
\mathbf{R}^d\rightarrow \mathbf{R}$ be a continuous function with
linear growth, that is there exist positive constants $K, D$, such
that $\forall( y, z)\in \mathbf{R}\times \mathbf{R}^d,\ |f(y,z)|\leq
K|y|+K|z|+D$. Then the sequence of functions
\begin{equation*}f_n(y,z)=\sup_{(y', z')\in Q\times Q^{d}}\{f(y',z')-n(|y-y'|+|z-z'|)\}
\end{equation*}
is well defined for $n \geq K$ and it satisfies\\
{\rm (i)}\hspace{4mm}  linear growth: $\forall ( y, z)\in
\mathbf{R}\times \mathbf{R}^d,
\ |f_n(y,z)|\leq K|y|+K|z|+D$;\\
{\rm (ii)}\hspace{3mm}  monotonicity in n: $\forall (y,z)\in
\mathbf{R}
\times \mathbf{R}^d,\ f_n(y,z)\searrow;$\\
{\rm (iii)} \hspace{2mm}Lipschitz condition: $\forall (y,z),(y',z')
\in \mathbf{R}\times \mathbf{R}^d$,
\begin{equation*}|f_n(y,z)-f_n(y',z')|\leq n(|y-y'|+|z-z'|);
\end{equation*}
{\rm (iv)}\hspace{1mm} strong convergence: if
$(y_n,z_n)\longrightarrow (y,z)$ then $f_n(y_n,z_n)\longrightarrow
f(y,z), n \rightarrow \infty. $
\end{lemma}

Consider, for fixed $(\omega ,t)$, the sequence $f_n(\omega,t,y,z)$
associated to $f$ by Lemma 4.2. Also consider
$F(\omega,t,y,z)=-|f(\omega,t,0,0)|-K|y|-K|z|.$ Then $f_n$ and
$F(\omega,t,y,z)$ are jointly measurable Lipschitz functions.
Given $\xi\in L^2(\Omega,{\mathcal{F}}_T,P)$, by Proposition 1,
there exist two pairs of processes $(Y^n,Z^n)$ and $(U,V)$ which are
the solutions of the following BDSDEs (4) and (5), respectively,
\begin{equation}
Y^n_t = \xi+\int_t^Tf_n(s,Y^n_s,Z^n_s){\rm d}s+\int_t^T
g(s,Y^n_s,Z^n_s){\rm d}B_s-\int_t^TZ^n_s{\rm d}W_s,
\end{equation}
\begin{equation}
U_t = \xi+\int_t^TF(s,U_s,V_s){\rm d}s+\int_t^Tg(s,U_s,V_s){\rm
d}B_s -\int_t^TV_s{\rm d}W_s.
\end{equation}
From the comparison theorem of \cite{SGL} and Lemma 4.2, we get
\begin{equation}
\forall n\geq m\geq K, \quad Y^m\geq Y^n\geq U, \quad {\rm d}t
\otimes {\rm d}P\mbox{-a.s.}
\end{equation}

\begin{lemma}
There exists a constant $A>0$ depending
only on $K$, $C$, $\alpha$, $T$ and $\xi$, such that
\begin{align*}
\forall n\geq K, \quad \|Y^n\|_{S^2}\leq A, \quad \|Z^n\|_{M^2}\leq
A; \quad \|U\|_{S^2}\leq A,\quad \|V\|_{M^2}\leq A.
\end{align*}
\end{lemma}

{\bf Proof.}\  First of all, we prove that $\|U\|_{S^2}$ and $\|V\|_{M^2}$
are all bounded. Clearly, from (6), there exists a constant $B$
depending only on $K$, $C$, $\alpha$, $T$ and $\xi$, such that
$$\left(\mathbf{E}\int_0^T|Y^n_s|^2{\rm d}s\right)^{1/2}\leq B, \quad \left(\mathbf{E}\int_0^T|U_s|^2{\rm d}s\right)^{1/2}\leq B,
\quad \|V\|_{M^2}\leq B.$$ Applying It\^o's formula to $|U_t|^2$, we
have
\begin{eqnarray}
\nonumber|U_t|^2 &=& |\xi|^2+2\int_t^TU_s\cdot F(s,U_s,V_s){\rm
d}s+2\int_t^TU_s\cdot g(s,U_s,V_s){\rm d}B_s\\
&&
-2\int_t^TU_s\cdot V_s{\rm d}W_s+\int_t^T|g(s,U_s,V_s)|^2{\rm d}s-\int_t^T|V_s|^2{\rm d}s.  %\eqno{\mbox{(4.5)}}%(13)
\end{eqnarray}
 From (H2), for all $\alpha<\alpha'<1$, there exists a
constant $C(\alpha')>0$ such that
\begin{eqnarray}
|g(t,u,v)|^2\leq C(\alpha')\left(|u|^2+|g(t,0,0)|^2\right)+\alpha'|v|^2.  %\eqno{\mbox{(4.6)}}           %(14)
\end{eqnarray}
From (7) and (8), it follows that
\begin{eqnarray*}
|U_t|^2+\int_t^T|V_s|^2{\rm d}s
& \leq & |\xi|^2+2K\int_t^T|U_s|(|f(s,0,0)|+|U_s|+|V_s|){\rm d}s\\
& & +C(\alpha')\int_t^T(|U_s|^2+|g(s,0,0)|^2){\rm d}s+\alpha'\int_t^T|V_s|^2{\rm d}s\\
& & +2\int_t^TU_s\cdot g(s,U_s,V_s){\rm d}B_s-2\int_t^TU_s\cdot V_s{\rm d}W_s\\
&\leq & |\xi|^2+K^2D(T-t)+C(\alpha')\int_t^T|g(s,0,0)|^2{\rm d}s\\
& &+\frac{1+\alpha'}{2}\int_t^T|V_s|^2{\rm d}s\\
& &+\left(1+2K+C(\alpha')+\frac{2K^2}{1-\alpha'}\right)\int_t^T|U_s|^2{\rm d}s\\
& &+2\int_t^TU_s\cdot g(s,U_s,V_s){\rm d}B_s-2\int_t^TU_s\cdot
V_s{\rm d}W_s.
\end{eqnarray*}
Taking supremum and expectation, by Young's inequality, we get
\begin{eqnarray}\nonumber
\|U\|_{S^2}^2+\frac{1-\alpha'}{2}\|V\|_{M^2}^2
& \leq &  \mathbf{E} \left(|\xi|^2+K^2TD+C(\alpha')\int_0^T|g(s,0,0)|^2{\rm d}s \right)\\
\nonumber
& &+\left(1+2K+C(\alpha')+\frac{2K^2}{1-\alpha'}\right) \mathbf{E}\int_0^T|U_s|^2{\rm d}s\\
\nonumber & & +2\mathbf{E}\sup_{0\leq t\leq T}\left|\int_t^TU_s\cdot
g(s,U_s,V_s){\rm d}B_s\right|\\
&&+2\mathbf{E}\sup_{0\leq t\leq T}\left|\int_t^TU_s\cdot V_s{\rm d}W_s\right|.%\quad \eqno{\mbox{(4.7)}} %(15)
\end{eqnarray}
By Burkholder-Davis-Gundy's inequality (cf. \cite{PP2}, \cite{SGL}), we deduce
\begin{eqnarray}\nonumber
&&\mathbf{E}\left(\sup_{0\leq t\leq T}\left|\int_t^TU_s\cdot
g(s,U_s,V_s){\rm d}B_s\right|\right)\\
\nonumber
& \leq &  c \mathbf{E}\left(\int_0^T|U_s|^2\cdot |g(s,U_s,V_s)|^2{\rm d}s\right)^{1/2}\\
\nonumber & \leq & c \mathbf{E}\left(\left(\sup_{0\leq t\leq T
}|U_t|^2\right)^{1/2}
\left(\int_0^T |g(s,U_s,V_s)|^2 {\rm d}s \right)^{1/2}\right)\\
\nonumber & \leq &2c^2C(\alpha')\mathbf{E} \left(\int_0^T|U_s|^2{\rm
d}s+\int_0^T |g(s,0,0)|^2 {\rm d}s\right )\\
&&+\frac18\|U\|_{S^2}^2+2c^2\alpha'\|V\|_{M^2}^2.
%\eqno{\mbox{(4.8)}}$$
\end{eqnarray} In the same way, we have
\begin{eqnarray}
\mathbf{E}\left(\sup_{0\leq t\leq T}\left|\int_t^TU_s\cdot V_s{\rm
d}W_s\right|\right)
\leq\frac18\|U\|_{S^2}^2+2c^2\|V\|_{M^2}^2.%\eqno{\mbox{(4.9)}} %(17)
\end{eqnarray}
From (10), (11) and (9), it follows that
\begin{eqnarray*}
&\quad & \|U\|_{S^2}^2+\frac{1-\alpha'}{2}\|V\|_{M^2}^2  \\
&\leq & 2\left(\mathbf{E}|\xi|^2+K^2TD+C(\alpha')(1+4c^2)\mathbf{E}\int_0^T|g(s,0,0)|^2{\rm d}s\right)\\
& & +2\left(1+2K+\frac{2K^2}{1-\alpha'}+C(\alpha')(1+4c^2)\right)\mathbf{E}\int_0^T|U_s|^2{\rm d}s\\
& &+8c^2(1+\alpha')\|V\|_{M^2}^2\\
&\leq & 2\left(\mathbf{E}|\xi|^2+K^2TD+C(\alpha')(1+4c^2)\mathbf{E}\int_0^T|g(s,0,0)|^2{\rm d}s\right)\\
& & +2\left(1+2K+\frac{2K^2}{1-\alpha'}+C(\alpha')(1+4c^2)+4c^2(1+\alpha')\right)B^2\\
& \mathbf{:=} & \frac{1-\alpha'}{2}(B')^2,
\end{eqnarray*}
that is
$$\|U\|_{S^2}\leq B',\quad \|V\|_{M^2}\leq B'.$$
From (6), it easily follows that
$$\quad \|Y^n\|_{S^2}\leq B'.$$
Next, we prove the boundedness  of $\|Z^n\|_{M^2}$. Applying It\^o's
formula to $|Y^n_t|^2$, it follows that
\begin{eqnarray*}
|Y^n_t|^2 &=& |\xi|^2+2\int_t^TY^n_s\cdot f_n(s,Y_s^n,Z_s^n){\rm d}s+2\int_t^TY^n_s\cdot g(s,Y^n_s,Z^n_s){\rm d}B_s\\
& &-2\int_t^TY^n_s\cdot Z^n_s{\rm
d}W_s+\int_t^T|g(s,Y^n_s,Z_s^n)|^2{\rm d}s-\int_t^T|Z_s^n|^2{\rm
d}s.
\end{eqnarray*}
Taking expectation, we have
\begin{eqnarray*}
\mathbf{E}(|Y^n_t|^2)+\mathbf{E}\int_t^T|Z^n_s|^2{\rm d}s
&=&\mathbf{E}|\xi|^2+2\mathbf{E}\int_t^TY^n_s\cdot
f_n(s,Y^n_s,Z^n_s){\rm d}s\\
&&+\mathbf{E}\int_t^T|g(s,Y^n_s,Z^n_s)|^2{\rm d}s.
\end{eqnarray*}
From  Young's inequality, it follows that
\begin{eqnarray*}
&&\mathbf{E}(|Y^n_t|^2)+\mathbf{E}\int_t^T|Z^n_s|^2{\rm d}s\\
&\leq& \mathbf{E}|\xi|^2+C'\mathbf{E}\int_t^T|Y^n_s|^2{\rm d}s
+\frac{1-\alpha'}{2}\mathbf{E}\int_t^T|Z^n_s|^2{\rm d}s\\
& + & K^2D(T-t)+C(\alpha')\mathbf{E}\int_t^T|g(s,0,0)|^2{\rm
d}s+\alpha' \mathbf{E}\int_t^T|Z^n_s|^2{\rm d}s,
\end{eqnarray*}
where $C'=1+2K+C(\alpha')+\dfrac{2K^2}{1-\alpha'}$, and we know
$0<\alpha'<1$ from (8). Then
\begin{eqnarray*}
\|Z^n\|_{M^2}^2&\leq&
\frac{2}{1-\alpha'}\left(C'T(B')^2+K^2DT+\mathbf{E}|\xi|^2
+C(\alpha')\mathbf{E}\int_0^T|g(s,0,0)|^2{\rm
d}s\right)\\&{\mathcal{:=}}&(A)^2.
\end{eqnarray*}
The proof is completed.
\quad$\Box$
\begin{lemma}
$\{(Y^n,Z^n)\}_{n=1}^{+\infty}$ converges in
$S^2([0,T];\mathbf{R})\times M^2(0,T;\mathbf{R}^d)$.
\end{lemma}
{\bf Proof.}\  Let $n_0\geq K$. Since $\{Y^n\}$ is decreasing
and bounded in $ S^2([0,T];\mathbf{R})$, we deduce from the
dominated convergence theorem that $Y^n$ converges in
$S^2([0,T];\mathbf{R})$. We shall denote by $Y$ the limit of
$\{Y^n\}$. Applying It\^o's formula to $|Y^n_t-Y^m_t|^2$, we get for
$n, m\geq n_0$,
\begin{eqnarray*}
& &\mathbf{E}(|Y^n_0-Y^m_0|^2)+\mathbf{E}\int_0^T|Z^n_s-Z^m_s|^2{\rm d}s \\
&=&2\mathbf{E}\int_0^T(Y^n_s-Y^m_s)(f_n(s,Y^n_s,Z^n_s)-f_m(s,Y^m_s,Z^m_s)){\rm d}s\\
& &+\mathbf{E}\int_0^T|g(s,Y^n_s,Z^n_s)-g(s,Y^m_s,Z^m_s)|^2{\rm d}s\\
&\leq& 2\left(\mathbf{E}\int_0^T |Y^n_s-Y^m_s|^2{\rm
d}s\right)^{\dfrac{1}{2}}
\left( \mathbf{E}\int_0^T|f_n(s,Y^n_s,Z^n_s)-f_m(s,Y^m_s,Z^m_s)|^2{\rm d}s\right)^{\dfrac{1}{2}}\\
& &+\mathbf{E}\int_0^T (C|Y^n_s-Y^m_s|^2+\alpha |Z^n_s-Z^m_s|^2
){\rm d}s.
\end{eqnarray*}
Since $f_n$ and $f_m$ are uniformly linear growth and
$\{(Y^n,Z^n)\}$ is bounded, similarly to Lemma 4.3, there exists a
constant $\bar{K}>0$ depending only on $K$, $C$, $\alpha$, $T$ and
$\xi$, such that
\begin{eqnarray*}
&&\mathbf{E}(|Y^n_0-Y^m_0|^2)+\mathbf{E}\int_0^T|Z^n_s-Z^m_s|^2{\rm
d}s\\ &\leq& \mathbf{E}\int_0^T( \bar{K}|Y^n_s-Y^m_s|^2+\alpha
|Z^n_s-Z^m_s|^2 ){\rm d}s.
\end{eqnarray*}
So
\begin{eqnarray*}
\|Z^n-Z^m\|_{M^2}^2\leq
\frac{\bar{K}T}{1-\alpha}\|Y^n-Y^m\|_{S^2}^2.
\end{eqnarray*}
Thus $\{Z^n\}$ is a Cauchy sequence in $M^2(0,T;\mathbf{R}^d)$, from
which the result follows. \quad$\Box$

{\bf Proof.}\  [Proof of Theorem \ref{result4}]
For all $n\geq n_0\geq K$ we have $Y^{n_0}\geq Y^n \geq U$, and
$\{Y^n\}$ converges in $S^2([0,T];\mathbf{R})$, ${\rm d}t\otimes
{\rm d}P$-a.s. to $Y\in S^2([0,T];\mathbf{R})$.

On the other hand, since $Z^n$ converges in $M^2(0,T;\mathbf{R}^d)$
to $Z$, we can assume, choosing a subsequence if needed, that
$Z^n\rightarrow Z$ ${\rm d}t\otimes {\rm d}P$-a.s. and
$\bar{G}=\sup_{n}|Z^n|$ is ${\rm d}t \otimes {\rm d}P $ integrable.
Therefore, from (i) and  (iv) of Lemma 4.2 we get for almost all
$\omega$,

\begin{align*}
f_n(t,Y^n_t,Z^n_t)&\longrightarrow f(t,Y_t,Z_t),
\quad (n\rightarrow \infty) \quad {\rm d}t\mbox{-a.e.}\\
|f_n(t,Y^n_t,Z^n_t)|&\leq |f(t,0,0)|+K\sup_n|Y^n_t|+K\sup_n|Z^n_t|\\
&=|f(t,0,0)|+K\sup_n|Y^n_t|+K\overline{G}_t.
\end{align*}
Thus, for almost all $\omega$ and uniformly in $t$, it holds that
\begin{equation*}
\int_t^Tf_n(s,Y^n_s,Z^n_s){\rm d}s\longrightarrow
\int_t^Tf(s,Y_s,Z_s){\rm d}s,\quad (n\rightarrow \infty).
\end{equation*} From the
continuity properties of the stochastic integral, it follows that
\begin{equation*}
\sup_{0\leq t\leq T} \left|\int_t^TZ^n_s{\rm d}W_s-\int_t^TZ_s{\rm
d}W_s \right|
 \longrightarrow 0 \quad \mbox{in probability},
\end{equation*}
\begin{equation*}
\sup_{0\leq t\leq T} \left|\int_t^Tg(s,Y^n_s,Z^n_s) {\rm
d}B_s-\int_t^Tg(s,Y_s,Z_s){\rm d}B_s \right|
 \longrightarrow 0 \quad \mbox{in probability}.
\end{equation*}
Choosing, again, a subsequence, we can assume that the above
convergence is $P$-a.s. Finally,
\begin{equation*}
 \begin{split}
    |Y^n_t-Y^m_t|\leq & \int_t^T|f_n(s,Y^n_s,Z^n_s)-f_m(s,Y_s^m,Z_s^m)|{\rm d}s\\
    &  + \left|\int_t^Tg(s,Y^n_s,Z^n_s){\rm d}B_s-\int_t^Tg(s,Y^m_s,Z^m_s){\rm d}B_s\right|\\
    &  + \left|\int_t^TZ^n_s{\rm d}W_s-\int_t^TZ^m_s{\rm d}W_s\right|,
 \end{split}
\end{equation*}
and taking limits on $m$ and supremum over $t$, we get
\begin{equation*}
 \begin{split}
   \sup_{0\leq t\leq T}|Y^n_t-Y_t|\leq & \int_0^T|f_n(s,Y^n_s,Z^n_s)-f(s,Y_s,Z_s)|{\rm d}s\\
    & +\sup_{0\leq t\leq T}\left|\int_t^Tg(s,Y^n_s,Z^n_s){\rm d}B_s-\int_t^Tg(s,Y_s,Z_s){\rm d}B_s\right|\\
    & +\sup_{0\leq t\leq T}\left|\int_t^TZ^n_s{\rm
d}W_s-\int_t^TZ_s{\rm d}W_s\right|,\quad P\mbox{-a.s.}
 \end{split}
\end{equation*}
from which it follows that $Y^n$ converges uniformly in $t$ to $Y$
(in particular, $Y$ is a continuous process). Note that $\{Y^n\}$ is
monotone; therefore, we actually have the uniform convergence for
the entire sequence and not just for a subsequence. Taking limits in
BDSDE (4), we deduce that $(Y,Z)$ is a solution of BDSDE (1).

Let $(\widehat{Y},\widehat{Z})\in S^2([0,T];\mathbf{R})\times
M^2(0,T;\mathbf{R}^d)$ is any solution of BDSDE (1). From Theorem
3.1, we get that $Y^n\geq \widehat{Y}$, $\forall n\in N$ and
therefore $Y\geq \widehat{Y}$ proves that $Y$ is  the maximal
solution. \quad$\Box$

\section{Kneser-type theorem for BDSDEs}

In this section, we shall discuss an interesting question: How many
solutions does a one-dimensional BDSDE (1) satisfying (H2)-(H4)
have? This is a classical Kneser-type problem.
Under some appropriate conditions, we shall
prove a Kneser-type theorem for BDSDEs satisfying (H2)-(H4).

\begin{theorem} [Kneser-type theorem for BDSDEs]\label{result5}
 We assume {\rm (H2)-(H4)}. Let
$(\underline{Y}_{t},\underline{Z}_{t})\in S^2([0,T];$
$\mathbf{R})\times M^2(0,T;\mathbf{R}^d)$ and
$(\overline{Y}_{t},\overline{Z}_{t})\in S^2([0,T];\mathbf{R})\times
M^2(0,T;\mathbf{R}^d)$ be the minimal and maximal solutions of BDSDE
{\rm (1)} with the terminal condition $\xi \in
L^{2}(\Omega,{\mathcal F}_{T},P;\mathbf{R})$. We also assume there
exists a function $z=h(t,y,\tilde z):\Omega\times [0,T]\times
\mathbf{R}\times \mathbf{R}^l\to \mathbf{R}^d$, which is the inverse
function of $\ \tilde z=g(t,y,z)$ with respect to $z$, and
$h(t,y,\tilde z)$ satisfies {\rm (H2)} with respect to $t, y,\tilde
z$. Then for any $t_{0}\in [0,T]$ and $\eta \in
L^{2}(\Omega,{\mathcal F}_{t_{0}},P;\mathbf{R})$ such that
\begin{align*}
\underline{Y}_{t_{0}}\leq \eta \leq \overline{Y}_{t_{0}},\quad
P\mbox{-a.s.},
\end{align*}
\noindent there is at least one solution
$(Y_{t},Z_{t})_{t\in[0,T]}\in S^2([0,T];\mathbf{R})\times
M^2(0,T;\mathbf{R}^d)$ of BDSDE {\rm (1)} satisfying
\begin{align*}
{Y}_{t_{0}}=\eta, \quad P\mbox{-a.s.}.
\end{align*}

Moreover the set of solutions of BDSDE {\rm (1)} is closed in $
S^2([0,T];\mathbf{R})\times M^2(0,T$; $\mathbf{R}^d).$ That is, for
any sequence of solutions $(Y^n_{t},Z^n_{t})_{t\in[0,T]}$ of BDSDE
{\rm (1)} for $n=1,...,$ if $(Y^n,Z^n)\rightarrow (Y,Z)$ in $
S^2([0,T];\mathbf{R})\times M^2(0,T;\mathbf{R}^d)$ as $n\rightarrow
\infty,$ then $(Y,Z)$ is also a solution of BDSDE {\rm (1)}.
\end{theorem}

{\bf Proof.}\  Let $(Y_{t}^1,Z_{t}^1)_{t\in[0,t_{0}]}\in
S^2([0,t_{0}];\mathbf{R})\times M^2(0,t_{0};\mathbf{R}^d)$ be a
solution of the following BDSDE
\begin{equation*}
Y_{t}^1 = \eta +\int_{t}^{t_{0}}f( s,Y_{s}^1,Z_{s}^1){\rm
d}s+\int_{t}^{t_{0}}g( s,Y_{s}^1,Z_{s}^1){\rm
d}B_{s}-\int_{t}^{t_{0}}Z_{s}^1{\rm d}W_{s},\ t\in[0,t_{0}].
\end{equation*}

Consider the following forward doubly stochastic differential
equation
\begin{equation}
Y_{t}^2 = \eta - \int_{t_{0}}^{t}f( s,Y_{s}^2,\widetilde
Z_{s}^2){\rm d}s-\int_{t_{0}}^{t}\widetilde Z_{s}^2{\rm
d}B_{s}+\int_{t_{0}}^{t}h( s,Y_{s}^2,\widetilde Z_{s}^2){\rm
d}W_{s},\ t\in[t_{0},T].
\end{equation}
By the assumptions, it is not difficult to see that (12) is a
forward ``BDSDE" on $[t_0,T]$, which satisfies the conditions of
Theorem 4.1. Consequently, there is a solution $(Y_{t}^2,\widetilde
Z_{t}^2)_{t\in[t_{0},T]}\in S^2([t_{0},T];\mathbf{R})\times
M^2(t_{0},T;\mathbf{R}^l)$ satisfying (12). We can define a
${\mathcal F}_{t}$-measurable time
\begin{equation*}
\tau=\inf\{t \geq
t_{0},Y_{t}^2\not\in(\underline{Y}_{t},\overline{Y}_{t})\}.
\end{equation*}
By $\underline{Y}_{T}=\overline{Y}_{T}=\xi$, we know that $\tau \leq
T.$ Now we define on [0,T]
\begin{equation*}\label{Multi}
 \begin{split}
    (Y_{t},Z_{t}) =& {\bf 1}_{[0,t_{0})}(t)(Y_{t}^1,Z_{t}^1)+
    {\bf 1}_{[t_{0},\tau)}(t)(Y_{t}^2,h(
    t,Y_{t}^2,\widetilde Z_{t}^2))\\
    &   +{\bf 1}_{[\tau,T]}(t)(\overline{Y}_{t},\overline{Z}_{t})
    {\bf 1}_{\{Y^2_{\tau}=\overline{Y}_{\tau}\}}
    +{\bf 1}_{[\tau,T]}(t)(\underline{Y}_{t},\underline{Z}_{t})
    {\bf 1}_{\{Y^2_{\tau}=\underline{Y}_{\tau}\}}.
 \end{split}
\end{equation*}
\noindent  It is easy to see that $(Y_{t},Z_{t})_{t\in[0,T]}\in
S^2([0,T];\mathbf{R})\times M^2(0,T;\mathbf{R}^d)$ is a solution of
BDSDE (1) with $Y_{T}=\xi$ and $Y_{t_{0}}=\eta.$

By the similar arguments in Theorem 4.1, it is easy to check the
closedness of the set of solutions of BDSDE (1) by the continuity of
$f$ with respect to $y$ and $z$.
\quad$\Box$

\begin{corollary} Indeed, in the case when the solution
of BDSDE {\rm (1)} is not unique, the cardinality of the set of the
associated solutions is at least of continuum since we can take
$$\eta=\alpha \underline{Y}_{t_0}+(1-\alpha)\overline{Y}_{t_0}$$ for
each $\alpha \in[0,1]$. Thus Theorem 5.1 is a Kneser-type theorem
for BDSDEs.
\end{corollary}

\section*{Acknowledgments} We would like to thank the anonymous referees
 for helpful comments and suggestions which led to much improvement of the earlier version of this paper.

% You may incorporate your references as follows in your main tex file.
% Using BibTex is not recommended but can be handled.

\end{document}